\documentclass[12pt]{article}
\usepackage{amsmath,amssymb}
\baselineskip 7pt \lineskip 7pt \numberwithin{equation}{section}

\def \Z{\hbox{$Z\hskip -5.2pt Z$}}

\def \C{\hbox{$C\hskip -5pt \vrule height 6pt depth 0pt \hskip 6pt$}}

\def\qed{\ \ \ifhmode\unskip\nobreak\fi\ifmmode\ifinner
         \else\hskip5pt\fi\fi
 \hbox{\hskip5pt\vrule width4pt height6pt depth1.5pt\hskip 1 pt}}

\def\l{\lambda}

\def\cl{\centerline}

\def\vs{\vspace*}

\def\C{\mathbb{C}}

\def\Z{\mathbb{Z}}

\begin{document}
\cl {{\large\bf Classification of irreducible weight modules with a
finite dimensional}}
 \cl {{\large\bf
 \vs{10pt}weight space over twisted Heisenberg-Virasoro algebra}\footnote{Supported by NSF grants 10471096, 10571120 of China
and ``One Hundred Talents Program'' from University of Science and
Technology of China}} \cl{Ran Shen$^*$,  \ Yucai
Su$^\dag$}\vs{4pt}\cl{\small\it $^*$Department of Mathematics,
Shanghai Jiao Tong University}\cl{\small\it Shanghai 200240,
China}\vs{4pt}\cl{\small\it \cl{\small\it $^\dag$Department of
Mathematics,
 University of Science and Technology of China}} \cl{\small\it Hefei 230026, China}
\vs{4pt}\cl{\small\it Email: ranshen@sjtu.edu.cn, \
ycsu@ustc.edu.cn}\vs{10pt}

{\small {\bf Abstract. } We show that the support of an irreducible
weight module over the twisted Heisenberg-Virasoro algebra, which
has an infinite dimensional weight space, coincides with the weight
lattice and that all nontrivial weight spaces of such a module are
infinite dimensional. As a corollary, we obtain that every
irreducible weight module over the twisted Heisenberg-Virasoro
algebra, having a nontrivial finite dimensional weight space, is a
Harish-Chandra module (and hence is either an irreducible highest or
lowest weight module or an irreducible module from the intermediate
series).
 \vs{2pt}\par {\bf Key
Words:} The twisted Heisenberg-Virasoro algebra, weight modules,
support \vs{2pt}\par{\it  Mathematics Subject Classification
(2000)}: 17B56; 17B68.}

\vs{10pt}\par \cl{\bf1. \ Introduction}
\setcounter{section}{1}\setcounter{theo}{0}
\par
\def\C{{\mathbb{C}}}
\def\Z{{\mathbb{Z}}}
\def\L{{\cal{L}}}
The {\it twisted Heisenberg-Virasoro algebra} is the Lie algebra
${\cal L}$ with  $\C$-basis $\{L_m,I_m,C,$ $C_I,C_{LI}\,|\, m \in
\mathbb{Z}\}$ subject to the following relations (e.g., [ACKP, B])
$$\begin{array}{lll} \!\!\!&\!\!\!&
[L_n,L_m]=(m-n)L_{n+m}+\delta_{n,-m}\frac{n^{3}-n}{12}C,
\\[7pt]\!\!\!&\!\!\!&
[L_n,I_m]=mI_{n+m}+\delta_{n,-m}(n^{2}+n)C_{LI}, \\[7pt]\!\!\!&\!\!\!&
[I_n,I_m]=n\delta_{n,-m}C_{I},\\[7pt] \!\!\!&\!\!\!& [{\mathcal
L},C]=[{\mathcal L},C_{LI}]=[{\mathcal L},C_{I}]=0.
\end{array}\eqno(1.1)$$
This Lie algebra is the universal central extension of the Lie
algebra of differential operators on a circle of order at most one,
which contains an infinite dimensional Heisenberg subalgebra and the
Virasoro subalgebra. The natural action of the Virasoro subalgebra
on the Heisenberg subalgebra is twisted with a 2-cocycle.

 Due to its
important applications in the representation theory of the toroidal
Lie algebras and some problems related to mathematical physics (for
instance, the even part of the well-known $N=2$ super-Virasoro
algebras is essentially the twisted Heisenberg-Virasoro algebra),
the structure and representation theories for the twisted
Heisenberg-Virasoro algebra has been well developed and has
attracted much attention in the literature (e.g., [ACKP, B, FO, JJ,
SJ, LZ]).

A {\it weight module} $M$ over ${\cal L}$ is an $\frak
H$-diagonalizable module, where $${\frak
H}:=\mbox{Span}_{\C}\{L_0,I_0,C, C_{LI},C_I\},$$ is the {\it Cartan
subalgebra} of $\L$. If, in addition, all weight spaces $M_\l$
(cf.~(1.2)) of a weight $\L$-module $M$ are finite dimensional, the
module is called a {\it Harish-Chandra module}. All irreducible
Harish-Chandra modules were classified in [LZ], and they are analogs
of that of the Virasoro algebra. Namely, they are exhausted by
irreducible highest weight modules, irreducible lowest weight
modules and irreducible modules from the so-called {\it intermediate
series} (see, e.g. [LZ]).

If $M$ is an irreducible weight $\L$-module, then $I_0,C,C_{LI},$
and $C_I$ must act as some complex numbers $h_I,c,c_I,$ and $
c_{LI}$ on $M$ respectively. Furthermore, $M$ has the weight space
decomposition
$$M=\bigoplus_{\lambda\in\C}M_\lambda, \mbox{ \ \ where }
M_\lambda=\{v\in M\,|\,L_0v=\lambda v\},\eqno(1.2)$$ where $M_\l$ is
called a {\it weight space with weight $\l$}. Denote
$${\rm Supp}(M):=\{\lambda\in\C\,|\,M_\lambda\neq 0\},$$
the set of all {\it weights $\l$} of $M$, called the {\it support}
of $M$. Obviously, if $M$ is an irreducible weight $\L$-module, then
there exists $\lambda\in\C$ such that ${\rm
Supp}(M)\subset\lambda+\Z$.

An irreducible weight module $M$ is called a {\it pointed module} if
there exists a weight $\l\in\C$ such that ${\rm dim\,}V_\l=1$. Xu
posted the following in [X]: \vs{4pt}\par{\bf Problem 1.1~}~{\it Is
any irreducible pointed module over the Virasoro algebra a
Harish-chandra module?}\vs{4pt}\par An irreducible weight module $M$
is called a {\it mixed  module} if there exist $\l\in\C$ and
$i\in\Z$ such that ${\rm dim\,}V_\l=\infty$ and ${\rm
dim\,}V_{\l+i}<\infty$. The following conjecture was posted in
[M]:\vs{4pt}\par {\bf Conjecture 1.2~}~{\it There are no irreducible
mixed module over the Virasoro algebra.}\vs{4pt}\par Mazorchuk and
Zhao [MZ] gave the positive answers to the above question and
conjecture.
 In this note, we also give the positive answers to the above
 question and conjecture for the twisted Heisenberg-Virasoro
 algebra. Our main result is the following:

\vs{4pt}\par {\bf Theorem 1.3~}~{\it Let M be an irreducible weight
$\L$-module. Assume that there exists $\lambda\in\C$ such that ${\rm
dim\,}M_\lambda=\infty$. Then ${\rm Supp}(M)=\lambda+\Z$, and for
every $k\in\Z$, we have ${\rm dim\,}M_{\lambda+k}=\infty$}.

Theorem 1.3 also implies the following classification of all
irreducible weight $\L$-modules which admit a nontrivial finite
dimensional weight space:

\vs{4pt}\par {\bf Corollary 1.4~}~{\it Let $M$ be an irreducible
weight $\L$-module. Assume that there exists $\lambda\in\C$ such
that $0<{\rm dim\,}M_\lambda<\infty$. Then $M$ is a Harish-Chandra
module. Consequently, $M$ is either an irreducible highest or lowest
weight module or an irreducible module from the intermidiate
series.} \vskip10pt
\par

\cl{\bf2. \ Proof of Theorem 1}
\setcounter{section}{2}\setcounter{theo}{0} \setcounter{equation}{0}
\vs{1pt}\par We first recall a main result about the weight
Virasoro-module in [MZ]:

\vs{4pt}\par {\bf Theorem 2.1~}~{\it Let $V$ be an irreducible
weight Virasoro-module. Assume that there exists $\lambda\in \C$,
such that ${\rm dim\,}V_\lambda=\infty$. Then ${\rm
Supp}(V)=\lambda+\Z$, and for every $k\in\Z$, we have ${\rm
dim\,}V_{\lambda+k}=\infty$.}

\vs{4pt}\par {\bf Lemma 2.2~}~{\it Assume that there exists
$\mu\in\C$ and a non-zero element $v\in M_\mu$, such that
$$I_1v=L_1v=L_2v=0\mbox{ \ \  or \ \ }I_{-1}v=L_{-1}v=L_{-2}v=0.$$ Then $M$ is a
Harish-Chandra module.}

\vs{4pt}\par{\it Proof.~}~Indeed, under these conditions, $v$ is
either a highest or a lowest weight vector, and hence, $M$ is a
Harish-Chandra module (see, e.g. [LZ]).$\hfill\Box$

Assume now that $M$ is an irreducible weight $\L$-module such that
there exists $\lambda\in\C$ satisfying ${\rm
dim\,}M_\lambda=\infty$.

\vs{4pt}\par {\bf Lemma 2.3~}~{\it There exists at most one $i\in\Z$
such that ${\rm dim\,}M_{\lambda+i}<\infty$.}

\vs{4pt}\par{\it Proof.~}~Assume that $${\rm
dim\,}M_{\lambda+i}<\infty\mbox{ \ \  and  \ }{\rm
dim\,}M_{\lambda+j}<\infty\mbox{ \ \ \ for some different \ \
}i,j\in\Z.$$ Without loss of generality, we may assume $i=1$ and
$j>1$. Set
$$
\begin{array}{ll}
V:=&\mbox{Ker}(I_1:M_\lambda\rightarrow M_{\lambda+1})\cap\mbox{Ker}
(L_1:M_\lambda\rightarrow
M_{\lambda+1})\cap\mbox{Ker}(I_j:M_\lambda\rightarrow
M_{\lambda+j})\\[7pt]&
\cap\,\mbox{Ker}(L_j:M_\lambda\rightarrow
M_{\lambda+j}),\end{array}$$ which is a subspace of $M_\l$. Since
$${\rm dim\,}M_\lambda=\infty,\ \ \ {\rm dim\,}M_{\lambda+1}<\infty\mbox{ \ \
and \ \ }{\rm dim\,}M_{\lambda+j}<\infty,$$ we have, ${\rm
dim\,}V=\infty.$ Since $$[L_1,L_k]=(k-1)L_{k+1}\neq 0\mbox{ \ \ and
\ }[I_1,L_l]=-I_{l+1}\neq 0\mbox{ \ \ for \ \ }k\in\Z,\ \ 1\le
l\in\Z,$$ we get
$$
\begin{array}{ll}
L_kV=0,& k=1,j,j+1,j+2,\cdots,\ \mbox{\ \ and \ \ }\\[7pt]
I_lV=0,& l=1,2,\cdots.\end{array}\eqno(2.1)$$If there would exist
$0\neq v\in V$ such that $L_2v=0$, then $I_1v=L_1v=L_2v=0$ and $M$
would be a Harish-Chandra module by Lemma 2.2. A contradiction.
Hence $L_2v\neq 0$ for all $v\in V$. In particular,
$${\rm dim\,}L_2V=\infty.$$ Since ${\rm dim\,}M_{\lambda+1}<\infty$,
and the actions of $I_{-1}$ and $L_{-1}$ on $L_2V$ map $L_2V$ (which
is an infinite dimensional subspace of $M_{\l+2}$) to $M_{\l+1}$
(which is finite dimensional), there exists $0\neq w\in L_2V$ such
that $I_{-1}w=L_{-1}w=0$. Let $w=L_2u$ for some $u\in V$. For all
$k\geq j$, using (2.1), we have
$$L_kw=L_kL_2u=L_2L_ku+(2-k)L_{k+2}u=0+0=0.$$ Hence $L_kw=0$ for
all $k=1,j,j+1,j+2,\cdots.$ Since
$$[L_{-1},L_l]=(l+1)L_{l-1}\neq 0\ \mbox{ \ and \ } [I_{-1},I_l]=I_{l-1}\neq 0\ \mbox{\ for
all \ }l>1,$$ we get inductively $L_kw=I_kw=0$ for all
$k=1,2,\cdots.$ Hence $M$ is a Harish-Chandra module by Lemma 2.2. A
contradiction. The lemma follows.$\hfill\Box$

Because of Lemma 2.3, we can now fix the following notation: $M$ is
an irreducible weight $\L$-module, $\mu\in\C$ is such that ${\rm
dim\,}M_\mu<\infty$ and ${\rm dim\,}M_{\mu+i}=\infty$ for every
$i\in\Z\setminus\{0\}$.

\vs{4pt}\par {\bf Lemma 2.4~}~{\it Let $0\neq v\in M_{\mu-1}$ be
such that $I_1v=L_1v=0$. Then

{\rm(1)} $(L_1^3-6L_2L_1+6L_3)L_2v=0,$

{\rm(2)} $I_kv=0,$

{\rm(3)} $I_kL_2v=0,$

\noindent where $k=1,2,\cdots.$}

\vs{4pt}\par{\it Proof.~}~(1) is the conclusion of Lemma 4 in [MZ].

(2) Since $I_1v=L_1v=0$ and $[I_l,L_1]=-lI_{l+1}\neq 0$ for $l\geq
1$, we inductively get (2).

(3) follows from $I_kL_2=L_2I_k-kI_{k+2}$ and (2).$\hfill\Box$

\vs{4pt}\par{\it Proof of Theorem 1.3.~}~Let
$$V:=\mbox{Ker}\{L_1:M_{\mu-1}\rightarrow M_\mu\}\cap
\mbox{Ker}\{I_1:M_{\mu-1}\rightarrow M_\mu\}\subset M_{\mu-1}.$$
Since ${\rm dim\,}M_{\mu-1}=\infty$ and ${\rm dim\,}M_\mu<\infty$,
we have ${\rm dim\,}V=\infty$. For any $v\in V$, consider the
element $L_2v$. By Lemma 2.2, $L_2v=0$ would imply that $M$ is a
Harish-Chandra module, a contradiction. Hence $L_2v\neq 0$, in
particular, ${\rm dim\,}L_2V=\infty$.

Since the actions of $I_{-1}$ and $L_{-1}$ on $L_2V$ map $L_2V$
(which is an infinite dimensional subspace of $M_{\mu+1}$) to
$M_\mu$ (which is finite dimensional), there exists $w=L_2v\in L_2V$
for some $v\in V$, such that $w\neq 0$ and $I_{-1}w=L_{-1}w=0$.
Using $[L_{-1},I_{-k+1}]=(1-k)I_{-k}$, we obtain $$I_{-k}w=0\mbox{ \
\ \ for \ \ } k=1,2,\cdots.\eqno(2.2)$$ From Lemma 2.4 (1) and (3),
we have $(L_1^3-6L_2L_1+6L_3)w=0,$ and $$I_kw=0,\
k=1,2,\cdots.\eqno(2.3)$$ Therefore,
$$I_0w=[L_{-1},I_1]w=0\mbox{ \ \ \ and \ \ }
C_Iw=[I_1,I_{-1}]w=0.\eqno(2.4)$$Using the commutation relations
(1.1), we have
$$I_{-1}^3(L_1^3-6L_2L_1+6L_3)=-48C_{LI}^3\mbox{ \ \ mod }\biggl(\bigoplus_{k\in\Z}
U(\L)I_k\biggr),$$ this implies $-48c_{LI}^3w=0$. Therefore,
$c_{LI}=0$. This together with (2.2)--(2.4) means that
$I_k,C_I,C_{LI}$ act trivially on $M$ for all $k\in\Z$, and so $M$
is simply a module over the Virasoro algebra. Thus, Theorem 1.3
follows from Theorem 2.1. $\hfill\Box$

\vs{4pt}\par{\it Proof of Corollary 1.4.~}~Assume that $M$ is not a
Harish-Chandra module. Then there should exists $i\in\Z$ such that
${\rm dim\,}M_{\lambda+i}=\infty$. In this case, Theorem 1.3 implies
${\rm dim\,}M_\lambda=\infty$, a contradiction. Hence $M$ is a
Harish-Chandra module, and the rest of the statement follows from
[LZ].$\hfill\Box$

\par
 \vs{10pt}
\par
\cl{\bf REFERENCES} \vs{-1pt}\small \lineskip=4pt
\begin{itemize}\parskip-5pt
\item[{[ACKP]}] E. Arbarello, C. De Concini, V.G. Kac, C. Procesi,
Moduli spaces of curves and representation theory, {\it Comm.
Math. Phys.}, {\bf117}(1988), 1-36.

\item[{[B]}] Y. Billig, Respresentations of the twisted
Heisenberg-Virasoro algebra at level zero, {\it Canad. Math.
Bulletin}, {\bf46}(2003), 529-537.

\item[{[FO]}] M.A. Fabbri, F. Okoh, Representations of
Virasoro-Heisenberg algebras and Virasoro-toroidal algebras, {\it
Canad. J. Math.}, {\bf51}(1999), no.3, 523-545.

\item[{[JJ]}] Q. Jiang, C. Jiang, Representations of the twisted
Heisenberg-Virasoro algebra and the full toroidal Lie algebras, {\it
Algebra Colloq.}, accepted.



\item[{[LZ]}] R. Lu, K. Zhao, Classification of irreducible weight
modules over the twisted Heisenberg-Virasoro algebra, preprint
(arXiv:math.ST/ 0510194)

\item[{[M]}] V. Mazorchuk,
On simple mixed modules over the Virasoro algebra, {\it Mat. Stud.}
{\bf22}(2004), 121-128.

\item[{[MZ]}] V. Mazorchuk, K. Zhao, Classification of simple
weight Virasoro modules with a finite dimensional weight space,
preprint (arXiv:math.RT/ 0507195).

\item[{[SJ]}] R. Shen, C. Jiang, Derivation algebra and
automorphism group of the twisted Heisenberg-Virasoro algebra,
Comm. Algebra, accepted.

\item[{[X]}] X. Xu, Pointed representations of Virasoro algebra,
A Chinese summary appears in Acta Math. Sinica 40 (1997) no. 3, 479.
{\it Acta Math. Sinica (N.S.)} {\bf13}(1997), 161-168.
\end{itemize}
\end{document}